\documentclass[preprint,12pt]{elsarticle}
\usepackage[utf8]{inputenc}

\usepackage{amsmath}
\usepackage{amssymb}
\usepackage{mathtools}
\usepackage{amsthm}
\usepackage{graphicx}
\usepackage{url}
\usepackage{multirow}
\usepackage{xcolor}

\mathtoolsset{showonlyrefs}
\allowdisplaybreaks[1]

\journal{}

\newtheorem{algorithm}{Algorithm}
\newtheorem{lemma}{Lemma}
\newtheorem{theorem}{Theorem}
\newtheorem{example}{Example}
\newtheorem{remark}{Remark}
\newtheorem{problem}{Problem}
\newcommand{\tg}[1]{#1}

\begin{document}

\maketitle

\begin{frontmatter}

\title{An efficient estimation of nested expectations without conditional sampling}
\author[label1]{Tomohiko Hironaka\corref{cor1}}
\ead{hironaka-tomohiko@g.ecc.u-tokyo.ac.jp}
\affiliation[label1]{organization={School of Engineering, The University of Tokyo},
            addressline={7-3-1 Hongo},
            city={Bunkyo-ku},
            postcode={113-8656},
            state={Tokyo},
            country={Japan}}
\author[label1]{Takashi Goda}
\ead{goda@frcer.t.u-tokyo.ac.jp}
\cortext[cor1]{Corresponding author}

\begin{abstract}
Estimating nested expectations is an important task in computational mathematics and statistics. In this paper we propose a new Monte Carlo method using post-stratification to estimate nested expectations efficiently without taking samples of the inner random variable from the conditional distribution given the outer random variable. This property provides the advantage over many existing methods that it enables us to estimate nested expectations only with a dataset on the pair of the inner and outer variables drawn from the joint distribution. We show an upper bound on the mean squared error of the proposed method under some assumptions. Numerical experiments are conducted to compare our proposed method with several existing methods (nested Monte Carlo method, multilevel Monte Carlo method, and regression-based method), and we see that our proposed method is superior to the compared methods in terms of efficiency and applicability.
\end{abstract}

\begin{highlights}
\item A new Monte Carlo estimator for nested expectations is proposed.
\item The proposed method has the advantage that it is free of inner conditional sampling.
\item An upper bound on the mean squared error is proven.
\item Numerical experiments confirm the superiority of the proposed methods over several existing methods.
\end{highlights}

\begin{keyword}
nested expectation\sep Monte Carlo method\sep expected value of sample information
\end{keyword}

\end{frontmatter}

\section{Introduction}\label{sec:intro}

Motivated by several different applications, we study estimating nested expectations which are defined as follows.
Let $X = (X_1, \cdots, X_J)\in\mathbb{R}^J$ and $Y = (Y_1, \cdots, Y_K)\in\mathbb{R}^K$ be possibly dependent random variables following the joint probability density $\rho(X,Y)$.
For a function $f: \mathbb{R}^J\to \mathbb{R}$, the nested expectation is defined by
\begin{equation}\label{eq:nested_expectation}
        I = \mathbb{E}_{\rho(Y)} f\left(\mathbb{E}_{\rho(X|Y)}X\right).
\end{equation}
Here we emphasize that the outer expectation is taken with respect to the marginal distribution of $Y$, while the inner expectation is with respect to the conditional distribution of $X$ given $Y$. Throughout this paper, we simply write $\rho(X)$ (resp.\ $\rho(Y)$) to denote the marginal probability density of $X$ (resp.\ $Y$), and also write $\rho(X|Y)$ (resp.\ $\rho(Y|X)$) to denote the conditional probability density of $X$ given $Y$ (resp.\ $Y$ given $X$).

The motivating examples are as follows:
\begin{example}[expected information gain]\label{exm:eig}
The concept of Bayesian experimental design aims to construct an optimal experimental design under which making the observation $Y$ maximizes the expected information gain (EIG) on the input random variable $\theta$ \cite{lindley1956measure,chaloner1995bayesian}.
Here the EIG denotes the expected amount of reduction in the Shannon information entropy and is given by
    \begin{align}
        & \mathbb{E}_{\rho(Y)} \mathbb{E}_{\rho(\theta|Y)}\log \rho(\theta|Y)-\mathbb{E}_{\rho(\theta)}\log \rho(\theta) \\
        & = \mathbb{E}_{\rho(\theta)}\mathbb{E}_{\rho(Y|\theta)}\log \rho(Y|\theta) - \mathbb{E}_{\rho(Y)}\log \left(\mathbb{E}_{\rho(\theta)}\rho(Y|\theta)\right),
    \end{align}
where the equality follows from Bayes' theorem. A nested expectation appears in the second term on the right-hand side.
\end{example}

\begin{example}[expected value of sample information]\label{exm:evsi}
Let $D$ be a finite set of possible medical treatments. As the outcome and cost of each treatment $d\in D$ is uncertain, we model its net benefit as a function of the input random variable $\theta$, denoted by $\mathrm{NB}_d$, where $\theta$ includes, for instance, the probability of side effect and the cost of treatment.

In the context of medical decision making, we want to know whether it is worth conducting a clinical trial or medical research to reduce the uncertainty of $\theta$ \cite{welton2012evidence}. Denoting the observation from a clinical trial or medical research by $Y$, the expected value of sample information (EVSI) measures the average gain in the net benefit from making the observation $Y$ and is given by
    \begin{equation}\label{eq:evsi}
        \mathbb{E}_{\rho(Y)}
            \max_{d\in D} \mathbb{E}_{\rho(\theta|Y)}\mathrm{NB}_d(\theta)
        - \max_{d\in D} \mathbb{E}_{\rho(\theta)}\mathrm{NB}_d(\theta),
    \end{equation}
where the first term represents the average net benefit when choosing the optimal treatment depending on the observation $Y$, and the second term does the net benefit without making the observation. Here the first term is exactly a nested expectation as given in \eqref{eq:nested_expectation}.
\end{example}

The nested Monte Carlo (NMC) method is probably the most straightforward approach to estimate nested expectations. The idea is quite simple: approximating the inner and outer expectations by the standard Monte Carlo methods, respectively. To be more precise, for positive integers $N_p$ and  $N_q$, the NMC estimator is given by
\begin{equation}\label{eq:nmc_estimator}
        \frac{1}{N_p}\sum_{p=1}^{N_p}f\left(
            \frac{1}{N_q}\sum_{q=1}^{N_q}X^{(p, q)}
        \right),
\end{equation}
where $Y^{(1)}, \ldots, Y^{(N_p)}$ denote the i.i.d.\ samples drawn from $\rho(Y)$, $X^{(p, 1)}, \cdots, X^{(p, N_q)}$ denote the i.i.d.\ samples drawn from $\rho(X|Y=Y^{(p)})$ for each $p=1,\ldots,N_p$, and the inner sum over $q$ is taken element-wise. However, it has been known that a large computational cost is necessary for the NMC method to estimate nested expectations with high accuracy \cite{rainforth2018nesting}. Moreover, the NMC method has a disadvantage in terms of applicability, since it requires generating the i.i.d.\ samples from $\rho(X|Y)$, which is often quite hard in applications \cite{strong2015estimating,goda2020multilevel,hironaka2020multilevel}.

A typical situation in estimating nested expectations is, instead, that we can generate i.i.d.\ samples from $\rho(X)$ and $\rho(Y|X)$, or, those from $\rho(X,Y)$. One way to tackle this issue is to use a Markov chain sampler directly for each $\rho(X|Y=Y^{(p)})$ in \eqref{eq:nmc_estimator}. Although the resulting estimator might be consistent, it is quite hard to obtain a non-asymptotic upper bound on the mean squared error and to choose an optimal allocation for $N_p$ and $N_q$ with the total cost fixed. Another way is to rewrite \eqref{eq:nested_expectation} into
\[ \mathbb{E}_{\rho(Y)} f\left(\frac{\mathbb{E}_{\rho(X)}X\rho(Y|X)}{\mathbb{E}_{\rho(X)}\rho(Y|X)}\right),\]
by using Bayes' theorem. The problem is then that, as a function of $X$, the likelihood $\rho(Y|X)$ is typically concentrated around a small region of the support of $X$, so that we need some sophisticated techniques like importance sampling to reliably estimate the ratio of two expectations inside.

There are some recent works to overcome these problems of the NMC estimator. We refer to \cite{goda2020multilevel,beck2018fast} for an efficient estimation of the EIG, and to \cite{strong2015estimating,hironaka2020multilevel,menzies2016efficient,jalal2018gaussian,heath2018efficient,heath2020calculating} for an efficient estimation of the EVSI. Estimating nested expectations is also an important task in portfolio risk measurement \cite{duffie2010dynamic, gordy2010nested}, and the relevant literature can be found in \cite{broadie2015risk,hong2017kernel,giles2019multilevel}

In this paper, we propose a new Monte Carlo estimator for nested expectations based on post-stratification, which does not require sampling from the conditional distribution \tg{$\rho(X|Y)$}. In fact, our proposed estimator only requires a dataset on the pair $(X,Y)$ drawn from the joint distribution $\rho(X,Y)$, and avoids any need to use a Markov chain sampler and importance sampling. Moreover, our proposed estimator is experimentally more efficient than several existing methods for numerical examples with small $K$ (the dimension of $Y$). The rest of this paper is organized as follows. In Section \ref{sec:method} we introduce our proposed estimator, and then give a theoretical analysis to derive an upper bound on the mean squared error of our proposed estimator in Section \ref{sec:theoretical}. Numerical experiments in Section \ref{sec:numerical} compares our proposed method with several existing methods (NMC method, multilevel Monte Carlo method \cite{goda2020multilevel,hironaka2020multilevel}, and regression-based method \cite{strong2015estimating}). We conclude this paper with discussion in Section~\ref{sec:discussion}.

\section{Method}\label{sec:method}

In this section, we first provide a new Monte Carlo estimator for nested expectations using post-stratification. Subsequently we introduce a variant of our proposed method using a linear regression.
\begin{algorithm}[proposed method]\label{alg:proposed}
    For $N = m ^ {2K}$ with $m \in\mathbb{N}$, let $
        \{(X^{(n)},Y^{(n)})\}_{n=1,\ldots,N}
    $ be a set of $N$ i.i.d.\ samples from the joint density $\rho(X,Y)$.

    \tg{
    For $k\in \{1,\ldots, K+1\}$, define one bijective map $t_k(\cdot, \cdot): \{1,\ldots,m^{k-1}\}\times \{1,\ldots,m^{2K-k+1}\}\to \{1,\ldots,N\}$ by
    \[ t_k(u,v)\coloneqq (u-1)m^{2K-k+1}+v.\]

    For $k\in \{0,\ldots, K\}$, define another bijective map $s_k: \{1,\ldots,N\}\to \{1,\ldots,N\}$ recursively such that $s_0(i)=i$, and for $k>0$, it satisfies
    \begin{align}
        \left\{ s_{k}(t_k(u, v))\right\}_{v\in \{1,\ldots,m^{2K-k+1}\}} & = \left\{ s_{k-1}(t_k(u, v))\right\}_{v\in \{1,\ldots,m^{2K-k+1}\}}
    \end{align}
    and
    \begin{align}
        Y_k^{(s_k(t_k(u, 1)))} \leq &\dots \leq Y_k^{(s_k(t_k(u, m^{2K-k+1})))}
    \end{align}
    for all $u\in \{1,\ldots,m^{k-1}\}.$ And let
    \begin{align}
        t(u, v) = s_K(t_{K+1}(u, v)),
    \end{align}
    for $u,v\in \{1,\ldots,m^{K}\}$.
    }

    Then our proposed estimator is given by
    \begin{equation}
        \widehat{I}
    \coloneqq
        \frac{1}{\sqrt{N}}
        \sum_{p=1}^{\sqrt{N}}
        f\left(
            \frac{1}{\sqrt{N}}
            \sum_{q=1}^{\sqrt{N}}
            X^{(\tg{t(p, q)})}
        \right),
    \end{equation}
    where the inner sum over $q$ is taken element-wise.
\end{algorithm}

\tg{
\begin{remark}\label{rem:set_equality}
 For each coordinate $k\in \{1,\ldots,K\}$, the following property holds. For any
    $p_1\in\{1, \dots, m^{k-1}\}$,
    $p_2\in\{1, \dots, m\}$, and
    $p_3\in\{1, \dots, m^{K-k}\}$, the sample set
    \[ \{ Y_k^{(t((p_1-1)m^{K-k+1}+(p_2-1)m^{K-k}+p_3, v))} \}_{v\in \{1,\ldots,\sqrt{N}\}}\]
    is a set of $\sqrt{N}$ points drawn from the \emph{parent} sample set
    \[ \{ Y_k^{(s_k(t_{k+1}((p_1-1)m+p_2, v)))}\}_{v\in \{1,\ldots,m^{2K-k}\}}\]
    consisting of $m^{2K-k}=\sqrt{N}m^{K-k}$ points, which coincides exactly with the set
    \[ \{ Y_k^{(s_k(t_{k}(p_1, (p_2-1)m^{2K-k} + v)))}\}_{v\in \{1,\ldots,m^{2K-k}\}}.\]
    This property plays a crucial role in the error analysis in the next section.
\end{remark}
}

The form of our proposed estimator apparently looks similar to the NMC estimator \eqref{eq:nmc_estimator} with $N_p=N_q=\sqrt{N}$. However, as we have mentioned already, the crucial difference is that our proposed estimator is free of the inner conditional sampling from the distribution $\rho(X|Y)$. We approximate the inner expectation $\mathbb{E}_{\rho(X|Y)}X$ by the average
\[ \tilde{X}^{(p)}=\frac{1}{\sqrt{N}}\sum_{q=1}^{\sqrt{N}} X^{(\tg{t(p, q)})},\]
based only on the set of $N$ i.i.d.\ samples from the joint distribution $\rho(X,Y)$, applying post-stratification to which is important in that the samples $X^{(\tg{t(p, 1)})},\ldots,X^{(\tg{t(p,\sqrt{N})})}$ are used as a substitute of the samples from the exact conditional distribution $\rho(X|Y)$. We note that, although each $X^{(\tg{n})}$ can be regarded as a sample from $\rho(X|Y=Y^{(\tg{n})})$, using only one inner sample is not enough to approximate the exact inner expectation. This is why we divide the set of $N$ samples into the $\sqrt{N}$ sets of $\sqrt{N}$ samples using post-stratification and use each subset to compute $\tilde{X}^{(p)}$.

Generating the i.i.d.\ samples from $\rho(X|Y)$ is often quite hard in some applications \cite{strong2015estimating,goda2020multilevel,hironaka2020multilevel}. Therefore, by construction, our proposed estimator has the advantage over many existing methods including the NMC estimator in terms of applicability. In the literature, there exist some other methods to estimate a special class of nested expectations, which do not rely on sampling from conditional distribution, see for instance \cite{strong2015estimating,hong2017kernel}. A regression-based method from \cite{strong2015estimating} is compared with our proposed estimator in Section~\ref{sec:numerical}.

We end this section with introducing a variant of Algorithm~\ref{alg:proposed}.
\begin{algorithm}[proposed method with regression]\label{alg:proposed_reg}
    For $N = m ^ {2K}$ with $m \in\mathbb{N}$, let $\{(X^{(n)},Y^{(n)})\}_{n=1,\ldots,N}$ be a set of $N$ i.i.d.\ samples from the joint density $\rho(X,Y)$.
    \tg{
    Let $s_k, t_k$ and $t$ be as defined in Algorithm~\ref{alg:proposed}.
    }
    Then our proposed estimator with a linear regression is given by
    \begin{equation}
        \widehat{I}_{\mathrm{rg}}
        \coloneqq \frac{1}{\sqrt{N}}\sum_{p=1}^{\sqrt{N}}
            \frac{1}{\sqrt{N}}\sum_{q=1}^{\sqrt{N}}
                f\left(\tilde{X}^{(\tg{t(p, q)})}\right),
    \end{equation}
    where $\tilde{X}^{(\tg{t(p, q)})}=(\tilde{X}^{(\tg{t(p, q)})}_1,\ldots,\tilde{X}^{(\tg{t(p, q)})}_J)$ is computed by
    \begin{equation}
        \tilde{X}_j^{(\tg{t(p, q)})} = \hat{c}^{\tg{(j, p)}}_0
            + \hat{c}^{\tg{(j, p)}}_1 Y_1^{(\tg{t(p, q)})}
            + \cdots
            + \hat{c}^{\tg{(j, p)}}_K Y_K^{(\tg{t(p, q)})}
    \end{equation}
    for each $j\in\{1, \cdots, J\}$ with
    \begin{equation}
        \hat{c}^{\tg{(j, p)}} = ( \hat{c}^{\tg{(j, p)}}_0,\hat{c}^{\tg{(j, p)}}_1,\ldots,\hat{c}^{\tg{(j, p)}}_K)^{\top}
    \end{equation}
    being a vector of coefficients which satisfies
    \begin{equation}
        M_{\tg{p}}^{\top}M_{\tg{p}}\hat{c}^{\tg{(j, p)}} = M_{\tg{p}}^{\top}v_{\tg{j, p}}
    \end{equation}
    with
    \begin{equation}
        M_{\tg{p}} = \begin{pmatrix}
            1      & Y_1^{(\tg{t(p, 1)})}        & \cdots & Y_k^{(\tg{t(p, 1)})}        & \cdots & Y_K^{(\tg{t(p, 1)})}       \\
            \vdots & \vdots              & \ddots &                     &        & \vdots             \\
            1      & Y_1^{(\tg{t(p, q)})}        &        & Y_k^{(\tg{t(p, q)})}        &        & Y_K^{(\tg{t(p, q)})}       \\
            \vdots & \vdots              &        &                     & \ddots & \vdots             \\
            1      & Y_1^{(\tg{t(p, \sqrt{N})})} & \cdots & Y_k^{(\tg{t(p, \sqrt{N})})} & \cdots & Y_K^{(\tg{t(p, \sqrt{N})})}
        \end{pmatrix}, \quad
        v_{\tg{j, p}} = \begin{pmatrix}
            X_j^{(\tg{t(p, 1)})} \\
            X_j^{(\tg{t(p, 2)})} \\
            \vdots \\
            X_j^{(\tg{t(p, \sqrt{N})})}
        \end{pmatrix}.
    \end{equation}
\end{algorithm}
Our motivation behind introducing this variant is that, instead of using the samples $X^{(\tg{t(p,1)})},\ldots,X^{(\tg{t(p,\sqrt{N})})}$ directly, we perturb them to approximate the inner conditional expectation $\mathbb{E}_{\rho(X|Y=Y^{(\tg{t(p,q)})})}X$ simply by one (perturbed) sample $\tilde{X}^{(\tg{t(p,q)})}$. Then the resulting estimator is given in a non-nested form.

\tg{
Both for Algorithms~\ref{alg:proposed} and \ref{alg:proposed_reg}, computing $s_k,t_k$ and $t$ requires $O(KN\log N)$ steps. In this paper, however, we assume that sampling $(X, Y)$ is the dominant part in computation, so that this cost for post-stratification is considered negligible and the number of total samples is used as a subjective measure of computational cost when comparing the performance of different estimation methods.
}

\section{Theoretical Analysis}\label{sec:theoretical}

In this section we show an upper bound on the mean squared error (MSE) of our proposed method. Here we show the result only for Algorithm~\ref{alg:proposed}. It is left open for future research whether a similar result on the MSE also holds for Algorithm~\ref{alg:proposed_reg}.  \tg{Throughout this section, the variance of a random vector $Z$ means the trace of the covariance matrix of $Z$ and we denote it by $\mathbb{V}_{\rho(Z)}Z$. This notation applies to the case of conditional variance with an appropriate subscript representing the underlying conditional distribution.}

\begin{theorem}[A bound of the MSE]\label{thm:proposed}
    Let $m$, $N$, $X^{(\tg{t(p, q)})}$, and $Y^{(\tg{t(p, q)})}$
    be given as in Algorithm~\ref{alg:proposed}.
    Let $F_Y$ be the cumulative distribution function of $Y$, and $F_{Y_k}$ be the marginal cumulative distribution function for $k\in\{1, 2, \dots, K\}$. Assume that every $F_{Y_k}$ is continuous, that $Y_j$ and $Y_k$ are mutually independent for all $1\leq j < k \leq K$, that there exist $\alpha, \beta>0$ such that
    \begin{align}\label{eq:assump1}
        \left|f(X_1) - f(X_2)\right| \leq \alpha\left\|X_1 - X_2\right\|_2
    \end{align}
    and
    \begin{align}\label{eq:assump2}
        \left\|\mathbb{E}_{\rho(X|Y_1)}X - \mathbb{E}_{\rho(X|Y_2)}X\right\|_2 \leq \beta\left|F_Y(Y_1) - F_Y(Y_2)\right|
    \end{align}
    hold for any $X_1, X_2\in\mathbb{R}^J$ and $Y_1, Y_2\in\mathbb{R}^K$, and that
    \begin{align}
        \mathbb{E}_{\rho(Y)}\mathbb{V}_{\rho(X|Y)}X < \infty.
    \end{align}
    Then we have
    \begin{align}
        \mathbb{E}(I-\widehat{I})^2 \leq 2\alpha^2\beta^2\frac{K^2}{N^{1/K}}
            + 2\alpha^2\left(\frac{\beta^2K}{3}+\mathbb{E}_{\rho(Y)}\mathbb{V}_{\rho(X|Y)}X\right)\frac{1}{N^{1/2}}.
    \end{align}
\end{theorem}

Before proving Theorem~\ref{thm:proposed}, we prove Lemma~\ref{lem:order}.

\begin{lemma}\label{lem:order}
    Let $m$, $N$, $X^{(\tg{t(p, q)})}$, and $Y^{(\tg{t(p, q)})}$
    be given as in Algorithm~\ref{alg:proposed}, and let $F_Y$ denote the cumulative distribution function of $Y$. Under the same assumptions considered in Theorem~\ref{thm:proposed}, we have
    \begin{align}\label{eq:lem1}
            \frac{1}{\sqrt{N}}\sum_{p=1}^{\sqrt{N}}
            \frac{1}{N}\sum_{q,q'=1}^{\sqrt{N}}
            \left(
                \mathbb{E}F_Y(Y^{(\tg{t(p, q)})}) - \mathbb{E}F_Y(Y^{(\tg{t(p, q')})})
            \right)^2
            \leq \frac{K^2}{N^{1/K}},
    \end{align}
    and
    \begin{align}\label{eq:lem2}
            \frac{1}{\sqrt{N}}\sum_{p=1}^{\sqrt{N}}
            \frac{1}{\sqrt{N}}\sum_{q=1}^{\sqrt{N}}
            \mathbb{V}F_Y(Y^{(\tg{t(p, q)})})
            \leq \frac{K}{6N^{1/2}}.
    \end{align}
\end{lemma}

\begin{proof}
    Let us recall the following two well-known facts as preparation. One fact is on order statistics. For a positive integer $n$, let $U_1,\ldots,U_n$ be i.i.d.\ random variables following the uniform distribution $U[0,1]$, and let $U_{(1,n)},\ldots,U_{(n,n)}$ be the corresponding order statistics. Then, for any $1\leq r\leq n$, the $r$-th order statistics $U_{(r,n)}$ follows the beta distribution $\mathrm{Beta}(r,n-r+1)$, which ensures \begin{align}\label{eq:order_stat0}
        \mathbb{V}U_{(r,n)}
            = \frac{r(n-r+1)}{(n+1)(n+2)}.
    \end{align}
    Moreover, for any $1\leq r<s\leq n$, the difference $U_{(s,n)} - U_{(r,n)}$ follows the beta distribution $\mathrm{Beta}(s-r,n-(s-r)+1)$ and so it holds that
    \begin{align}\label{eq:order_stat}
        \mathbb{E}\left(U_{(s,n)} - U_{(r,n)}\right)^2
            = \frac{(s-r)(s-r+1)}{(n+1)(n+2)}.
    \end{align}

    The other fact is on an elementary inequality. For a positive integer $n$ and any $\{a_1,\ldots,a_n\},\{b_1,\ldots,b_n\}\subset [0,1]$, it holds that
    \begin{align}\label{eq:elem_inequ}
    \left| \prod_{i=1}^{n}a_i-\prod_{i=1}^{n}b_i\right|\leq \sum_{i=1}^{n}\left| a_i-b_i\right|.
    \end{align}
    This inequality can be proven by an induction on $n$. As the case $n=1$ trivially holds, let us assume that the inequality holds for $n=m$ and prove the case $n=m+1$. In fact,
    \begin{align}
        \left| \prod_{i=1}^{m+1}a_i-\prod_{i=1}^{m+1}b_i\right| & = \left| (a_{m+1}-b_{m+1})\prod_{i=1}^{m}a_i+b_{m+1}\left(\prod_{i=1}^{m}a_i- \prod_{i=1}^{m}b_i\right)\right| \\
        & \leq \left| a_{m+1}-b_{m+1}\right|\prod_{i=1}^{m}a_i+b_{m+1}\left|\prod_{i=1}^{m}a_i- \prod_{i=1}^{m}b_i\right| \\
        & \leq \left| a_{m+1}-b_{m+1}\right|+\sum_{i=1}^{m}\left| a_i-b_i\right|=\sum_{i=1}^{m+1}\left| a_i-b_i\right|,
    \end{align}
    \tg{where the first inequality follows from the triangle inequality and the second inequality follows from the induction assumption and the trivial bound $\prod_{i=1}^{m}a_i\leq 1$.} This completes the proof of \eqref{eq:elem_inequ}.

    Now let us prove \eqref{eq:lem1}. To do so, let us consider a single coordinate $Y_k$ for $k\in \{1,\ldots,K\}$ first. In what follows, we denote by $F_{Y_k}$ the marginal cumulative distribution function of $Y_k$ as stated in Theorem~\ref{thm:proposed}. Here, since every one of the outer variables $Y_k$ is assumed independent from all the other variables $Y\setminus Y_k$, we can see that, for any $p_1\in\{1, \dots, m^{k-1}\}$, the sorted random variables $F_k(Y_k^{(\tg{s_k(t_k(p_1, 1))})}), \dots, F_k(Y_k^{(\tg{s_k(t_k(p_1, m^{2K-k+1}))})})$ correspond to the order statistics $U_{(1,m^{2K-k+1})},\ldots,U_{(m^{2K-k+1},m^{2K-k+1})}$, respectively.

    Therefore, for $k\in \{1,\ldots,K\}$, we have
        \begin{align}
        &\quad
            \frac{1}{\sqrt{N}}\sum_{p=1}^{\sqrt{N}}
            \frac{1}{N}\sum_{q,q'=1}^{\sqrt{N}}
            \mathbb{E}
            \left(
                F_{Y_k}(Y_k^{(\tg{t(p, q)})})
                - F_{Y_k}(Y_k^{(\tg{t(p, q')})})
            \right)^2 \\
        &=
            \frac{1}{m^{k-1}}\sum_{p_1=1}^{m^{k-1}}
            \frac{1}{m}\sum_{p_2=1}^{m}
            \frac{1}{m^{K-k}}\sum_{p_3=1}^{m^{K-k}}
            \frac{1}{N}\sum_{q,q'=1}^{\sqrt{N}}
            \\
        &\quad
            \mathbb{E}
            \left(
                F_{Y_k}(Y_k^{(\tg{t((p_1-1)m^{K-k+1}+(p_2-1)m^{K-k}+p_3, q)})})
                - F_{Y_k}(Y_k^{(\tg{t((p_1-1)m^{K-k+1}+(p_2-1)m^{K-k}+p_3, q')})})
            \right)^2 \\
        &=
            \frac{1}{m^{k-1}}\sum_{p_1=1}^{m^{k-1}}
            \frac{1}{m}\sum_{p_2=1}^{m}
            \frac{1}{m^{2(2K-k)}}\sum_{q,q'=1}^{m^{2K-k}}
            \\
        &\quad
            \mathbb{E}
            \left(
                F_{Y_k}(Y_k^{(\tg{s_k(t_k(p_1, (p_2-1)m^{2K-k}+q))})})
                - F_{Y_k}(Y_k^{(\tg{s_k(t_k(p_1, (p_2-1)m^{2K-k}+q'))})})
            \right)^2 \\
        &=
            \frac{1}{m^{k-1}}\sum_{p_1=1}^{m^{k-1}}
            \frac{1}{m}\sum_{p_2=1}^{m}
            \frac{1}{m^{2(2K-k)}}\sum_{q,q'=1}^{m^{2K-k}}
            \mathbb{E}
            \left(
                U_{((p_2-1)m^{2K-k}+q, m^{2K-k+1})}
                - U_{((p_2-1)m^{2K-k}+q', m^{2K-k+1})}
            \right)^2 \\
        &=
            \frac{1}{m^{k-1}}\sum_{p_1=1}^{m^{k-1}}
            \frac{1}{m}\sum_{p_2=1}^{m}
            \frac{1}{m^{2(2K-k)}}\sum_{q,q'=1}^{m^{2K-k}}
            \frac{|q-q'|(|q-q'|+1)}{(m^{2K-k+1}+1)(m^{2K-k+1}+2)} \\
        &= \frac{(m^{2K-k}-1)m^{2K-k}(m^{2K-k}+1)(m^{2K-k}+2)}
            {6m^{2(2K-k)}(m^{2K-k+1}+1)(m^{2K-k+1}+2)} \leq \frac{1}{m^2}=\frac{1}{N^{1/K}},\label{eq:each_dimension}
    \end{align}
    \tg{where we used the identity
    \begin{align}
    \{1, \ldots, m^{K}\} & = \{(p_1-1)m^{K-k+1}+(p_2-1)m^{K-k}+p_3 \\
    & \qquad \mid p_1\in \{1,\ldots,m^{k-1}\}, p_2\in \{1,\ldots,m\}, p_3\in \{1,\ldots,m^{K-k}\} \},
    \end{align}
    to obtain the first equality, the second equality follows from Remark~\ref{rem:set_equality}, the third equality follows from the argument on the order statistics made in the above paragraph, and the fourth equality follows from \eqref{eq:order_stat}.}

    Now we get
    \begin{align}
        &\quad
            \frac{1}{\sqrt{N}}\sum_{p=1}^{\sqrt{N}}
            \frac{1}{N}\sum_{q,q'=1}^{\sqrt{N}}
            \left(
                \mathbb{E}F_Y\left(Y^{(\tg{t(p, q)})}\right)
                - \mathbb{E}F_Y\left(Y^{(\tg{t(p, q')})}\right)
            \right)^2 \\
        &=
            \frac{1}{\sqrt{N}}\sum_{p=1}^{\sqrt{N}}
            \frac{1}{N}\sum_{q,q'=1}^{\sqrt{N}}
            \left(
                \prod_{k=1}^{K}\mathbb{E}F_{Y_k}(Y_k^{(\tg{t(p, q)})})
                - \prod_{k=1}^{K}\mathbb{E}F_{Y_k}(Y_k^{(\tg{t(p, q')})})
            \right)^2 \\
        &\leq
            \frac{1}{\sqrt{N}}\sum_{p=1}^{\sqrt{N}}
            \frac{1}{N}\sum_{q,q'=1}^{\sqrt{N}}
            \left(\sum_{k=1}^{K}\left|
                \mathbb{E}F_{Y_k}(Y_k^{(\tg{t(p, q)})})
                - \mathbb{E}F_{Y_k}(Y_k^{(\tg{t(p, q')})})\right|
            \right)^2 \\
        &\leq
            K\sum_{k=1}^{K}\frac{1}{\sqrt{N}}\sum_{p=1}^{\sqrt{N}}
            \frac{1}{N}\sum_{q,q'=1}^{\sqrt{N}}
            \left|
                \mathbb{E}F_{Y_k}(Y_k^{(\tg{t(p, q)})})
                - \mathbb{E}F_{Y_k}(Y_k^{(\tg{t(p, q')})})\right|^2 \\
        &\leq \frac{K^2}{N^{1/K}},
    \end{align}
    \tg{where the first equality follows from the assumption that every $Y_k$ is mutually independent, the first and second inequalities follow from \eqref{eq:elem_inequ} and Jensen's inequality, respectively, and the last bound follows from \eqref{eq:each_dimension}.} This completes the proof of \eqref{eq:lem1}.

Let us move on to proving \eqref{eq:lem2}. By a
reasoning similar to what is used to prove \eqref{eq:lem1}, we have
    \begin{align}
        &\quad
            \frac{1}{\sqrt{N}}\sum_{p=1}^{\sqrt{N}}
            \frac{1}{\sqrt{N}}\sum_{q=1}^{\sqrt{N}}
            \mathbb{V}F_{Y_k}(Y_k^{(\tg{t(p, q)})}) \\
        &=  \frac{1}{m^{k-1}}\sum_{p_1=1}^{m^{k-1}}
            \frac{1}{m}\sum_{p_2=1}^{m}
            \frac{1}{m^{2K-k}}\sum_{q=1}^{m^{2K-k}}
            \mathbb{V}U_{((p_2-1)m^{2K-k}+q, m^{2K-k+1})}\\
        &=  \frac{1}{m^{2K-k+1}}\sum_{q=1}^{m^{2K-k+1}}
            \mathbb{V}U_{(q, m^{2K-k+1})}\\
        &= \frac{1}{m^{2K-k+1}}\sum_{q=1}^{m^{2K-k+1}}
            \frac{q(m^{2K-k+1}-q+1)}{(m^{2K-k+1}+1)^2(m^{2K-k+1}+2)} \\
        &= \frac{1}{6(m^{2K-k+1}+1)}, \label{eq:sum_variance}
        \end{align}
        \tg{where the first equality follows from Remark~\ref{rem:set_equality} and the previous argument on the order statistics and we used \eqref{eq:order_stat0} to obtain the third equality.}

    Therefore, it holds that
    \begin{align}
        & \quad
            \frac{1}{\sqrt{N}}\sum_{p=1}^{\sqrt{N}}
            \frac{1}{\sqrt{N}}\sum_{q=1}^{\sqrt{N}}
            \mathbb{V}F_Y(Y^{(\tg{t(p, q)})}) \\
        & = \frac{1}{\sqrt{N}}\sum_{p=1}^{\sqrt{N}}
            \frac{1}{\sqrt{N}}\sum_{q=1}^{\sqrt{N}}
            \mathbb{E}\left(\prod_{k=1}^{K}F_{Y_k}(Y_k^{(\tg{t(p, q)})})-\prod_{k=1}^{K}\mathbb{E}F_{Y_k}(Y_k^{(\tg{t(p, q)})})\right)^2\\
        & \leq \frac{1}{\sqrt{N}}\sum_{p=1}^{\sqrt{N}}
            \frac{1}{\sqrt{N}}\sum_{q=1}^{\sqrt{N}}
            \mathbb{E}\left(\sum_{k=1}^{K}\left|F_{Y_k}(Y_k^{(\tg{t(p, q)})})-\mathbb{E}F_{Y_k}(Y_k^{(\tg{t(p, q)})})\right|\right)^2\\
        &  \leq K\sum_{k=1}^{K}
            \frac{1}{\sqrt{N}}\sum_{p=1}^{\sqrt{N}}
            \frac{1}{\sqrt{N}}\sum_{q=1}^{\sqrt{N}}
            \mathbb{E}\left|F_{Y_k}(Y_k^{(\tg{t(p, q)})})-\mathbb{E}F_{Y_k}(Y_k^{(\tg{t(p, q)})})\right|^2 \\
        &   = K\sum_{k=1}^{K}
            \frac{1}{\sqrt{N}}\sum_{p=1}^{\sqrt{N}}
            \frac{1}{\sqrt{N}}\sum_{q=1}^{\sqrt{N}}
            \mathbb{V}F_{Y_k}(Y_k^{(\tg{t(p, q)})})\\
        &   = K\sum_{k=1}^{K}\frac{1}{6(m^{2K-k+1}+1)}
        \leq \frac{K(m^K-1)}{6m^{2K}(m-1)}\leq \frac{K}{6m^K}=\frac{K}{6N^{1/2}},
    \end{align}
    \tg{where the first equality comes from the assumption that every $Y_k$ is mutually independent, the first and second inequalities follow from \eqref{eq:elem_inequ} and Jensen's inequality, respectively, and we used \eqref{eq:sum_variance} in the third equality.} This completes the proof of \eqref{eq:lem2}.
\end{proof}

Now we are ready to show a proof of  Theorem~\ref{thm:proposed}.

\begin{proof}[Proof of Theorem~\ref{thm:proposed}]
Throughout this proof, we simply denote by $\mathbb{E}_{\bullet}$ the expectation with respect to the underlying probability measure of the corresponding random variable $\bullet$. First we obtain
    \begin{align}
        &\quad \mathbb{E}(I-\widehat{I})^2 \notag \\
        & = \mathbb{E}_{\{(X^{(\tg{t(p, q)})},Y^{(\tg{t(p, q)})})\}}
            \left(
                \mathbb{E}_Y f\left(\mathbb{E}_{X|Y}X\right)
                - \frac{1}{\sqrt{N}}\sum_{p=1}^{\sqrt{N}} f\left(
                    \frac{1}{\sqrt{N}}\sum_{q=1}^{\sqrt{N}}X^{(\tg{t(p, q)})}
                \right)
            \right)^2 \notag \\
        & = \mathbb{E}_{\{(X^{(\tg{t(p, q)})},Y^{(\tg{t(p, q)})})\}}
            \left(
                \mathbb{E}_Y f\left(\mathbb{E}_{X|Y}X\right)
                - \frac{1}{\sqrt{N}}\sum_{p=1}^{\sqrt{N}} f\left(
                    \frac{1}{\sqrt{N}}\sum_{q=1}^{\sqrt{N}}\mathbb{E}_{X|Y^{(\tg{t(p, q)})}}X
                \right) \right. \notag \\
        &\quad \quad \left. + \frac{1}{\sqrt{N}}\sum_{p=1}^{\sqrt{N}} f\left(
                    \frac{1}{\sqrt{N}}\sum_{q=1}^{\sqrt{N}}\mathbb{E}_{X|Y^{(\tg{t(p, q)})}}X
                \right)-\frac{1}{\sqrt{N}}\sum_{p=1}^{\sqrt{N}} f\left(
                    \frac{1}{\sqrt{N}}\sum_{q=1}^{\sqrt{N}}X^{(\tg{t(p, q)})}
                \right)
            \right)^2 \notag \\
        & \leq 2\mathbb{E}_{\{Y^{(\tg{t(p, q)})}\}}
            \left(
                \mathbb{E}_Y f\left(\mathbb{E}_{X|Y}X\right)
                - \frac{1}{\sqrt{N}}\sum_{p=1}^{\sqrt{N}} f\left(
                    \frac{1}{\sqrt{N}}\sum_{q=1}^{\sqrt{N}}\mathbb{E}_{X|Y^{(\tg{t(p, q)})}}X
                \right)
            \right)^2 \\
        &\quad + 2\mathbb{E}_{\{(X^{(\tg{t(p, q)})},Y^{(\tg{t(p, q)})})\}}
            \left(
                \frac{1}{\sqrt{N}}\sum_{p=1}^{\sqrt{N}} f\left(
                    \frac{1}{\sqrt{N}}\sum_{q=1}^{\sqrt{N}}\mathbb{E}_{X|Y^{(\tg{t(p, q)})}}X
                \right)
                - \frac{1}{\sqrt{N}}\sum_{p=1}^{\sqrt{N}} f\left(
                    \frac{1}{\sqrt{N}}\sum_{q=1}^{\sqrt{N}}X^{(\tg{t(p, q)})}
                \right)
            \right)^2, \label{eq:proof_thm_1}
        \end{align}
        \tg{where the first equality is given by the definition, to which we add and subtract a common term in the second equality, and the last bound follows from Jensen's inequality.}

        In what follows, we show an upper bound on each term on the right-most side of \eqref{eq:proof_thm_1}. Let us consider the first term. By denoting the i.i.d.\ copy of the outer random variables $\{Y^{(\tg{t(p,q)})}\}$ by $\{Y'^{(\tg{t(p,q)})}\}$, we have
        \begin{align}
        & \quad \mathbb{E}_{\{Y^{(\tg{t(p, q)})}\}}
            \left(
                \mathbb{E}_Y f\left(\mathbb{E}_{X|Y}X\right)
                - \frac{1}{\sqrt{N}}\sum_{p=1}^{\sqrt{N}} f\left(
                    \frac{1}{\sqrt{N}}\sum_{q=1}^{\sqrt{N}}\mathbb{E}_{X|Y^{(\tg{t(p, q)})}}X
                \right)
            \right)^2\\
        & = \mathbb{E}_{\{Y^{(\tg{t(p, q)})}\}}
            \left(
                \mathbb{E}_{\{Y'^{(\tg{t(p, q)})}\}} \frac{1}{\sqrt{N}}\sum_{p=1}^{\sqrt{N}}\frac{1}{\sqrt{N}}\sum_{q=1}^{\sqrt{N}}f\left(\mathbb{E}_{X|Y'^{(\tg{t(p, q)})}}X\right)
                - \frac{1}{\sqrt{N}}\sum_{p=1}^{\sqrt{N}} f\left(
                    \frac{1}{\sqrt{N}}\sum_{q=1}^{\sqrt{N}}\mathbb{E}_{X|Y^{(\tg{t(p, q)})}}X
                \right)
            \right)^2\\
        &\leq \mathbb{E}_{\{Y^{(\tg{t(p, q)})}\},\{Y'^{(\tg{t(p, q)})}\}}
            \frac{1}{\sqrt{N}}\sum_{p=1}^{\sqrt{N}}\frac{1}{\sqrt{N}}\sum_{q=1}^{\sqrt{N}}
            \left(
                f\left(\mathbb{E}_{X|Y'^{(\tg{t(p, q)})}}X\right)
                - f\left(
                    \frac{1}{\sqrt{N}}\sum_{q'=1}^{\sqrt{N}}\mathbb{E}_{X|Y^{(\tg{t(p, q')})}}X
                \right)
            \right)^2 \\
        &\leq \alpha^2\mathbb{E}_{\{Y^{(\tg{t(p, q)})}\},\{Y'^{(\tg{t(p, q)})}\}}
            \frac{1}{\sqrt{N}}\sum_{p=1}^{\sqrt{N}}\frac{1}{\sqrt{N}}\sum_{q=1}^{\sqrt{N}}
            \left\|
                \mathbb{E}_{X|Y'^{(\tg{t(p, q)})}}X
                - \frac{1}{\sqrt{N}}\sum_{q'=1}^{\sqrt{N}}\mathbb{E}_{X|Y^{(\tg{t(p, q')})}}X
            \right\|_2^2 \\
        &\leq \alpha^2\mathbb{E}_{\{Y^{(\tg{t(p, q)})}\},\{Y'^{(\tg{t(p, q)})}\}}
            \frac{1}{\sqrt{N}}\sum_{p=1}^{\sqrt{N}}\frac{1}{N}\sum_{q,q'=1}^{\sqrt{N}}
            \left\|
                \mathbb{E}_{X|Y'^{(\tg{t(p, q)})}}X - \mathbb{E}_{X|Y^{(\tg{t(p, q')})}}X
            \right\|_2^2 \\
        &\leq \alpha^2\beta^2
            \mathbb{E}_{\{Y^{(\tg{t(p, q)})}\},\{Y'^{(\tg{t(p, q)})}\}}
            \frac{1}{\sqrt{N}}\sum_{p=1}^{\sqrt{N}}
            \frac{1}{N}\sum_{q,q'=1}^{\sqrt{N}}
            \left(F_Y(Y'^{(\tg{t(p, q)})}) - F_Y(Y^{(\tg{t(p, q')})})\right)^2 \\
        &\leq \alpha^2\beta^2
            \frac{1}{\sqrt{N}}\sum_{p=1}^{\sqrt{N}}
            \frac{1}{N}\sum_{q,q'=1}^{\sqrt{N}}
            \left(\mathbb{E}F_Y(Y^{(\tg{t(p, q)})}) - \mathbb{E}F_Y(Y^{(\tg{t(p, q')})})\right)^2 + 2\alpha^2\beta^2
            \frac{1}{\sqrt{N}}
            \sum_{p=1}^{\sqrt{N}}
            \frac{1}{\sqrt{N}}\sum_{q=1}^{\sqrt{N}}
                \mathbb{V}F_Y(Y^{(\tg{t(p, q)})})\\
        &\leq \alpha^2\beta^2\frac{K^2}{N^{1/K}}
            + \alpha^2\beta^2\frac{K}{3N^{1/2}},
    \end{align}
    \tg{where Jensen's inequality was used in the first and third inequalities, and the second and fourth inequalities follow from the assumptions \eqref{eq:assump1} and \eqref{eq:assump2}, respectively, and the last bound follows from the results of Lemma~\ref{lem:order}.}

Now let us consider the second term on the right-most side of \eqref{eq:proof_thm_1}. Similarly to the above, we have
        \begin{align}
        &\quad  \mathbb{E}_{\{(X^{(\tg{t(p, q)})},Y^{(\tg{t(p, q)})})\}}
            \left(
                \frac{1}{\sqrt{N}}\sum_{p=1}^{\sqrt{N}} f\left(
                    \frac{1}{\sqrt{N}}\sum_{q=1}^{\sqrt{N}}\mathbb{E}_{X|Y^{(\tg{t(p, q)})}}X
                \right)
                - \frac{1}{\sqrt{N}}\sum_{p=1}^{\sqrt{N}} f\left(
                    \frac{1}{\sqrt{N}}\sum_{q=1}^{\sqrt{N}}X^{(\tg{t(p, q)})}
                \right)
            \right)^2 \\
        &\leq \mathbb{E}_{\{(X^{(\tg{t(p, q)})},Y^{(\tg{t(p, q)})})\}}
            \frac{1}{\sqrt{N}}\sum_{p=1}^{\sqrt{N}}
            \left(
                f\left(
                    \frac{1}{\sqrt{N}}\sum_{q=1}^{\sqrt{N}}\mathbb{E}_{X|Y^{(\tg{t(p, q)})}}X
                \right)
                - f\left(
                    \frac{1}{\sqrt{N}}\sum_{q=1}^{\sqrt{N}}X^{(\tg{t(p, q)})}
                \right)
            \right)^2 \\
        &\leq \alpha^2\mathbb{E}_{\{(X^{(\tg{t(p, q)})},Y^{(\tg{t(p, q)})})\}}
            \frac{1}{\sqrt{N}}\sum_{p=1}^{\sqrt{N}}
            \left\|
                \frac{1}{\sqrt{N}}\sum_{q=1}^{\sqrt{N}}\mathbb{E}_{X|Y^{(\tg{t(p, q)})}}X
                - \frac{1}{\sqrt{N}}\sum_{q=1}^{\sqrt{N}}X^{(\tg{t(p, q)})}
            \right\|_2^2 \\
        &= \tg{\alpha^2
            \mathbb{E}_{\{Y^{(\tg{t(p, q)})}\}}
            \frac{1}{\sqrt{N}}\sum_{p=1}^{\sqrt{N}}\frac{1}{N}\sum_{q=1}^{\sqrt{N}}
            \mathbb{V}_{X|Y^{(\tg{t(p, q)})}}X} \\
        &= \alpha^2\frac{1}{N^{1/2}}
            \mathbb{E}_{Y}\mathbb{V}_{X|Y}X,
        \end{align}
        \tg{where Jensen's inequality was used in the first inequality, and the second inequality follows from the assumption \eqref{eq:assump1}, and the first equality follows from the conditional independence between $X^{(n)}$ and $X^{(n')}$ for any $n\neq n'$ given $\{Y^{(n)}\}_{1\leq n\leq N}$, and the second equality follows from the fact that the sample set $\{Y^{(\tg{t(p, q)})}\}_{1\leq p, q \leq \sqrt{N}} = \{Y^{(n)}\}_{1\leq n\leq N}$ is the set of i.i.d.\ samples of $Y$.}

        Altogether we obtain the desired bound on the MSE as shown in Theorem~\ref{thm:proposed}.
\end{proof}

Our bound shows that the MSE decays at the rate of $N^{-1/\max(2,K)}$ at worst. It remains a challenge whether we can develop an algorithm without conditional sampling which achieves the $K$-independent rate of the MSE. Still it is interesting that our proposed estimator is comparable to the nested Monte Carlo estimator \eqref{eq:nmc_estimator} for $K\leq 2$, see \cite[Theorem~1]{rainforth2018nesting}, although the latter requires a necessity for generating i.i.d.\ samples of the inner variable from the conditional density $\rho(X|Y)$. Moreover, our theoretical analysis assumes the mutual independence between the individual outer variables. Thus relaxing the necessary conditions to get a similar result on the MSE is another interesting research question.

\section{Numerical Experiments}\label{sec:numerical}

Here we conduct numerical experiments to evaluate the performance of our proposed method and compare with some existing methods.

\subsection{Problem Settings}

As we have introduced in Section~\ref{sec:intro} as the motivating examples on nested expectations, we consider one toy problem related to EIG and two testcases to estimate EVSI, respectively. In what follows, we first show the model setting of each of the three problems.

\begin{problem}[toy problem for EIG from \cite{rainforth2018nesting}]\label{prm:1}
    Let $\theta,Y\in\mathbb{R}$ be independent random variables.
    We want to estimate the following nested expectation
    \begin{equation}
        I = \mathbb{E}_{\rho(Y)} \log\left(\mathbb{E}_{\rho(\theta)}
            \sqrt{\frac{2}{\pi}} \exp\left(-2(Y-\theta)^2\right)
        \right)
    \end{equation}
    where $Y\sim\mathrm{U}(-1,1)$ and $\theta\sim\mathrm{N}(0,1)$.
\end{problem}

\begin{problem}[simple testcase for EVSI]\label{prm:2}
    Let $\theta\in\mathbb{R}$ and $Y=(Y_1, Y_2, Y_3)\in\mathbb{R}^3$
    be random variables which follow the multivariate normal distribution:
    \begin{equation}
        \begin{pmatrix}
            \theta \\
            Y_1 \\
            Y_2 \\
            Y_3
        \end{pmatrix} \sim \mathrm{Norm}\left(
            \begin{pmatrix}
                0 \\
                0 \\
                0 \\
                0
            \end{pmatrix}, \begin{pmatrix}
                1 & \frac{1}{2} & \frac{1}{2} & \frac{1}{2} \\
                \frac{1}{2} & 1 & \frac{1}{2} & \frac{1}{2} \\
                \frac{1}{2} & \frac{1}{2} & 1 & \frac{1}{2} \\
                \frac{1}{2} & \frac{1}{2} & \frac{1}{2} & 1
            \end{pmatrix}
        \right)
    \end{equation}
    and let $D=\{1,2\}$. The net benefit functions $\mathrm{NB}_1$ and $\mathrm{NB}_2$ are given by
    \begin{equation}
        \mathrm{NB}_1(\theta) = \theta,\qquad \mathrm{NB}_2(\theta) = -\theta,
    \end{equation}
    respectively. We want to estimate the EVSI as defined in \eqref{eq:evsi}.
\end{problem}

\begin{problem}[medical decision problem for EVSI from \cite{hironaka2020multilevel}]\label{prm:3}
    Let $\theta$ be a vector consisting of 12 random variables, the distribution of each of which is described in Table~\ref{table:h2020}.
    \begin{table}[t]
        \caption{Model inputs involved in our experiments from \cite{hironaka2020multilevel}. Note that $\text{log-normal}(\mu,\Sigma)$ and $\text{logit-normal}(\mu,\Sigma)$ denote the log-normal and logit-normal distributions, respectively, with $\mu$ and $\Sigma$ being the mean vector and the covariance matrix of the corresponding normal distribution, respectively. The logit function $\text{logit}(p)$ is the transformation mapping probability $p$ to the logarithm of the odds $\log(p/(1-p))$.}\label{table:h2020}
        \small{ \begin{center}
        \begin{tabular}{|l|l|l|}
            \hline
            Description & Parameter & Distribution \\
            \hline \hline
            Lifetime remaining & $L$ & $N(30, 25)$ \\ \hline
            QALY after critical event, per year & $Q_E$ & $\text{logit-normal}\left(0.6, 1/36\right)$ \\ \hline
            QALY decrement due to side effects & $Q_{SE}$ & $N(0.7, 0.01)$ \\ \hline
            Cost of critical event & $C_E$ & $N(2\times 10^5, 10^8)$ \\ \hline
            Cost of side effect & $C_{SE}$ & $N(10^5, 10^8)$ \\ \hline
            Cost of treatment $d=1$ & $C_{T,1}$ & 0 (constant) \\ \hline
            Cost of treatments $d=2,3$ & $C_{T,d}$ & \multirow{2}{*}{$N\left(\begin{pmatrix} 1.5\times 10^4 \\ 2\times 10^4\end{pmatrix}, \begin{pmatrix}300 & 100 \\ 100 & 500\end{pmatrix}\right)$} \\
            & & \\ \hline
            Probability of critical event & $P_{E,1}$ & $\text{Beta}(15, 85)$ \\
            on treatment $d=1$ & &  \\ \hline
            Odds ratios of critical event & $OR_{E,d}$ & \multirow{2}{*}{$\text{log-normal}\left(\begin{pmatrix} -1.5 \\ -1.75\end{pmatrix}, \begin{pmatrix}0.11 & 0.02 \\ 0.02 & 0.06\end{pmatrix}\right)$} \\
            relative to treatment $d=1$ & &  \\
            $(P_{E,d}/(1-P_{E,d}))/(P_{E,1}/(1-P_{E,1}))$ & & \\ \hline
            Probability of critical event & $P_{E,d}$ & Derived from\\
            on treatments $d=2,3$ & & $P_{E,1}$ and $OR_{E,d}$\\ \hline
            Probability of side effect & $P_{SE,1}$ & 0 (constant) \\
            on treatment $d=1$ & & \\ \hline
            Probability of side effect & $P_{SE,d}$ & \multirow{2}{*}{$\text{logit-normal}\left(\begin{pmatrix} -1.4 \\ -1.1\end{pmatrix}, \begin{pmatrix}0.10 & 0.05 \\ 0.05 & 0.25\end{pmatrix}\right)$} \\
             on treatments $d=2,3$ & & \\ \hline
            Monetary value of 1 QALY & $\lambda$ & \$75,000 (constant) \\
            \hline
        \end{tabular} \end{center} }
    \end{table}
    Let $D = \{1, 2, 3\}$, and for each $d\in D$ we define the net benefit function as:
    \begin{align}
        \mathrm{NB}_d(\theta) &= P_{SE,d}P_{E,d}\left[\lambda\left(L\frac{1+Q_E}{2}-Q_{SE}\right)-(C_{SE}+C_E)\right] \\
            & \quad + P_{SE,d}(1-P_{E,d})\left[\lambda(L-Q_{SE})-C_{SE}\right] \\
            & \quad + (1-P_{SE,d})P_{E,d}\left[\lambda L\frac{1+Q_E}{2}-C_E\right] \\
            & \quad + (1-P_{SE,d})(1-P_{E,d})\lambda L-C_{T,d}.
    \end{align}
    Here we consider a small two-arm randomized controlled trial (RCT)
    comparing treatments $d=1$ and $d=3$ with $n_p=100$ patients,
    which informs three parameters $OR_{E,3}$, $C_{T,3}$ and $P_{SE,3}$.
    The sample information $Y$ is defined
    as the three-dimensional vector $(Y_1,Y_2,Y_3)$ which follows
    \begin{align}
        Y_1 &\sim N(\log(OR_{E,3}), 4/n_p), \\
        Y_2 &\sim N(C_{T,3}, 10^4/n_p), \\
        Y_3 &\sim \mathrm{Binomial}(n_p, P_{SE,3}),
    \end{align}
    respectively. Again, we want to estimate the EVSI as defined in \eqref{eq:evsi}.
\end{problem}

The first two models are so simple that the exact values of $I$ and EVSI can be calculated analytically, respectively. Therefore, these models can be used to verify that our proposed estimator is a consistent estimator for nested expectations, i.e., the estimate converges to a correct value as $N\to \infty$. On the other hand, the third model problem is more complicated and the analytical value of EVSI is not available. We use this model to test that our proposed estimator can work well even for complex models. Regarding the second problem, the outer variables are not mutually independent. \tg{Although our theoretical analysis does not deal with such situation, our proposed algorithm can be applied anyway.}

Note that, in the medical decision problem, the i.i.d.\ samples from the conditional distribution $\rho(\theta|Y)$ cannot be generated. Nevertheless, our proposed method can be applied, whereas some other methods, including the NMC estimator, cannot be directly applied without any help of Markov chain sampler or importance sampling.

In this paper, we compare our proposed method (without/with regression) with the existing three methods: the NMC method (for the first two problems), the multilevel Monte Carlo (MLMC) method \cite{hironaka2020multilevel} (for the third problem), and the regression-based method \cite{strong2015estimating}  (for all the three problems). In particular, we use a generalized additive model (GAM) for the regression-based method, so that we simply denote it by the GAM-based method; we confirmed beforehand that using a Gaussian process regression is too expensive to complete the experiments within a reasonable computational time. The source code used in our experiments is available at
\url{https://github.com/Goda-Research-Group/experiment-2021-09}.

\subsection{Results}

\begin{figure}
    \centering
    \includegraphics[width=\linewidth]{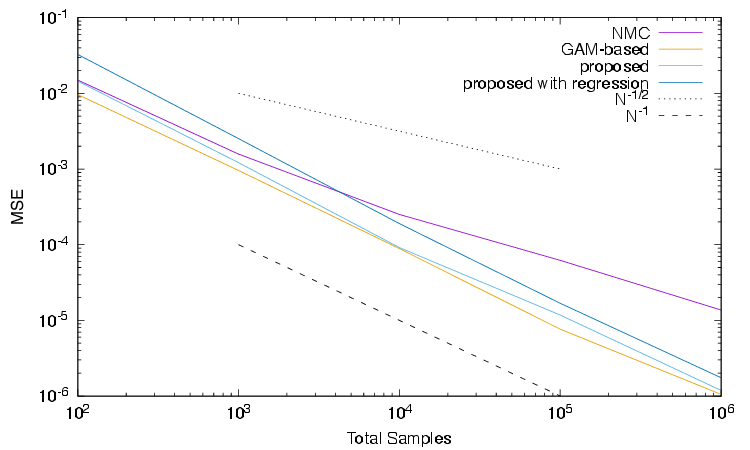}
    \caption{\tg{Comparison on MSE of the four methods (the NMC method, the GAM-based method, and our proposed methods without and with regression) for Problem~\ref{prm:1}. MSE is empirically estimated by averaging the results over 100 independent runs. The observed convergence rates are almost $O(N^{-1})$ for three methods excluding the NMC method and $O(N^{-1/2})$ for the NMC method.}}
    \label{fig:result_ne_testcase_comparison}
\end{figure}

For the first experiment to estimate $I$, we compare the four methods (the NMC method, the GAM-based method, and our proposed methods without and with regression) by the MSE against the total number of samples. The result is shown in Figure \ref{fig:result_ne_testcase_comparison}, where the vertical axis represents the MSE, and the horizontal axis does the total number of samples used. Here the MSE is empirically estimated by averaging the results over $100$ independent runs for each setting. In this simple test problem, \tg{the estimated slopes of the graphs (which show the convergence rates) are $-0.74806$ for NMC method, $-1.0204$ for proposed method, $-1.0733$ for proposed method with regression, and $-1.0733$ for GAM-based method.} We can see that the efficiencies of the GAM-based method and our proposed method without regression are quite comparable, and are better than the NMC method especially when the total number of samples is large. Our proposed method with regression is less efficient than the two methods, but the difference from them is not so large and becomes smaller as the number of samples increases.

\begin{figure}
    \centering
    \includegraphics[width=\linewidth]{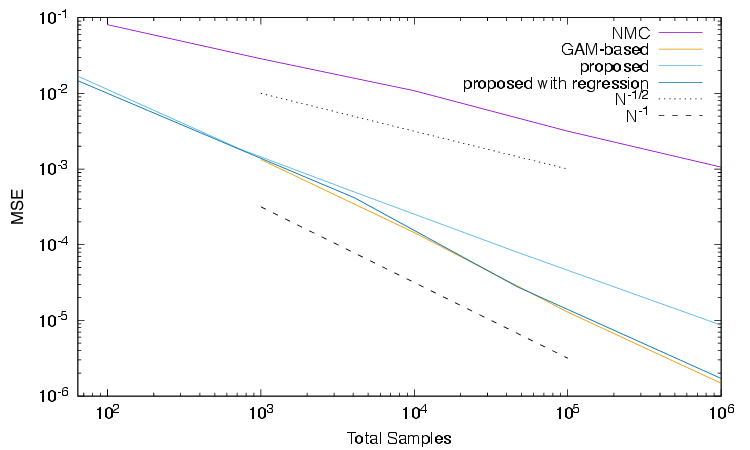}
    \caption{\tg{Comparison on MSE of the four methods (the NMC method, the GAM-based method, and our proposed methods without and with regression) for Problem~\ref{prm:2}. MSE is empirically estimated by averaging the results over 100 independent runs. The observed convergence rates are almost $O(N^{-1})$ for the GAM-based method and our proposed method with regression and $O(N^{-1/2})$ for the NMC method and our proposed method without regression.}}
    \label{fig:result_evsi_testcase_comparison}
\end{figure}

For the second experiment we compare the same four methods (the NMC method, the GAM-based method, our proposed methods without and with regression) as for the first experiment. The result is shown in Figure \ref{fig:result_evsi_testcase_comparison}. The meanings of axes are the same as Figure \ref{fig:result_ne_testcase_comparison}. \tg{The estimated slopes of the graphs are $-0.47191$ for NMC method, $-0.77749$ for proposed method, $-0.95066$ for proposed method with regression, and $-0.99162$ for GAM-based method.} Figure \ref{fig:result_evsi_testcase_comparison} clearly shows that the NMC method is inefficient compared to the other three methods. For this model, the GAM-based method and our proposed method with regression perform best among the four methods. The improvement over our proposed method without regression increases as the number of samples increases.

\begin{figure}
    \centering
    \includegraphics[width=\linewidth]{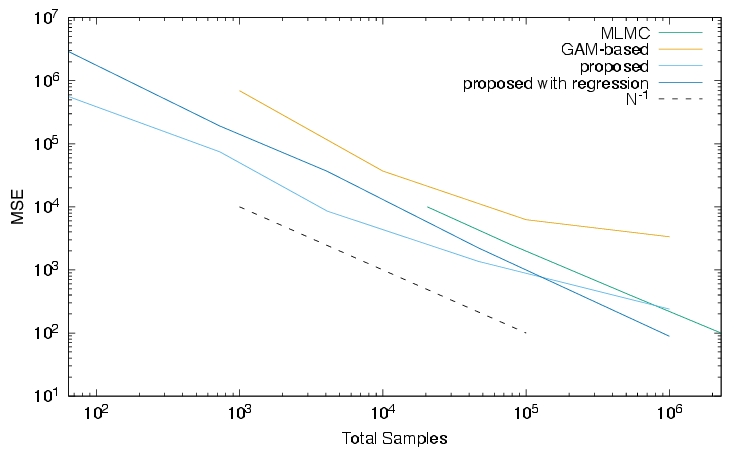}
    \caption{\tg{Comparison on MSE of the four methods (the MLMC method, the GAM-based method, and our proposed methods without and with regression) for Problem~\ref{prm:3}. MSE is empirically estimated by averaging the results over 100 independent runs. The GAM-based method does not converge to the correct value. The observed convergence rates are almost $O(N^{-1})$ for the MLMC method and our proposed method with regression.}}
    \label{fig:result_evsi_medical_comparison}
\end{figure}

Finally we compare the four methods (the MLMC method, the GAM-based method, our proposed methods without and with regression) for the last, medical decision problem. Although the exact value of EVSI is required for estimating the MSE, we cannot calculate EVSI analytically for this model. Therefore we use the value of EVSI estimated with high accuracy in \cite{hironaka2020multilevel} as a substitute, which is 1031. The result is shown in Figure \ref{fig:result_evsi_medical_comparison}. \tg{The estimated slopes of the graphs are $-0.82133$ for proposed method, $-1.0764$ for proposed method with regression, $-0.77180$ for GAM-based method, and $-0.97346$ for MLMC method.} For this model, our proposed method without and with regression provide the best results. When the number of samples is relatively small, our proposed method without regression performs better,
whereas, when the number of samples is large, using regression helps to improve the performance.

\section{Discussion}\label{sec:discussion}

We conclude this paper with some discussions on comparison of our proposed method with the other three methods.

Firstly, our proposed method outperforms the NMC method for the first two problems. Considering that the NMC method requires generating the i.i.d.\ samples from the conditional density $\rho(X|Y)$ whereas our proposed method does not, our proposed method not only is more efficient (at least for the case of small $K$) but also has a wider applicability.

Secondly, comparing our proposed method with the MLMC method for the third problem, the difference of the MSE is relatively small, so that our proposed method does not have a clear advantage in terms of efficiency. However, the MLMC method used to estimate the EVSI, developed in \cite{hironaka2020multilevel}, requires a good approximation of the conditional density $\rho(X|Y)$ as an importance distribution. Thus  it can be claimed that our proposed method has a wider applicability in this sense. Moreover, \tg{although the MLMC method is proven to achieve the MSE convergence rate $O(1/N)$ in theory,} it is not clear whether it does work well in practice when the total number of samples is small. This is because a certain amount of samples is required as an overhead to get a rough estimate on the variance of stochastic quantities. Therefore our proposed method can hold a superiority over the MLMC method either when we are in a pre-asymptotic region (i.e., the total number of samples is small) or a good importance distribution is not available.

Finally, let us look at the GAM-based method. We can see that it performs quite well for two simple test cases, but is inferior to the other methods for the third medical decision problem. GAM is a method of approximating a multi-variate function by a sum of uni-variate functions. In the context of nested expectations, $\mathbb{E}_{\rho(X|Y)}X$ is approximated by $\sum_{k=1}^{K}g_k(Y_k)$ for regressed functions $g_1,\ldots,g_K$. Although this assumption is satisfied for the first two testcases, it does not for the third medical decision problem. Therefore the GAM-based estimator even does not converge to the correct value. It is shown in the theoretical analysis (Theorem \ref{thm:proposed}) that our proposed method converges to the correct value even if such additivity of $f$ does not hold. This is a clear advantage of our proposed method. Since the our proposed method performs quite comparably with the GAM-based method even for the two simple test cases, we can say that our proposed method is superior to the GAM-based method.

\section*{Funding}
The work of the second author is supported by JSPS KAKENHI Grant Number 20K03744.

\bibliographystyle{elsarticle-num}
\bibliography{contents/91-references}

\begin{thebibliography}{10}
\expandafter\ifx\csname url\endcsname\relax
  \def\url#1{\texttt{#1}}\fi
\expandafter\ifx\csname urlprefix\endcsname\relax\def\urlprefix{URL }\fi
\expandafter\ifx\csname href\endcsname\relax
  \def\href#1#2{#2} \def\path#1{#1}\fi

\bibitem{lindley1956measure}
D.~V. Lindley, On a measure of the information provided by an experiment, The
  Annals of Mathematical Statistics (1956) 986--1005.

\bibitem{chaloner1995bayesian}
K.~Chaloner, I.~Verdinelli, Bayesian experimental design: A review, Statistical
  Science (1995) 273--304.

\bibitem{welton2012evidence}
N.~J. Welton, A.~J. Sutton, N.~Cooper, K.~R. Abrams, A.~Ades, Evidence
  synthesis for decision making in healthcare, Vol. 132, John Wiley \& Sons,
  2012.

\bibitem{rainforth2018nesting}
T.~Rainforth, R.~Cornish, H.~Yang, A.~Warrington, F.~Wood, On nesting {M}onte
  {C}arlo estimators, in: International Conference on Machine Learning, PMLR,
  2018, pp. 4267--4276.

\bibitem{strong2015estimating}
M.~Strong, J.~E. Oakley, A.~Brennan, P.~Breeze, Estimating the expected value
  of sample information using the probabilistic sensitivity analysis sample: a
  fast, nonparametric regression-based method, Medical Decision Making 35~(5)
  (2015) 570--583.

\bibitem{goda2020multilevel}
T.~Goda, T.~Hironaka, T.~Iwamoto, Multilevel {M}onte {C}arlo estimation of
  expected information gains, Stochastic Analysis and Applications 38~(4)
  (2020) 581--600.

\bibitem{hironaka2020multilevel}
T.~Hironaka, M.~B. Giles, T.~Goda, H.~Thom, Multilevel {M}onte {C}arlo
  estimation of the expected value of sample information, SIAM/ASA Journal on
  Uncertainty Quantification 8~(3) (2020) 1236--1259.

\bibitem{beck2018fast}
J.~Beck, B.~M. Dia, L.~F. Espath, Q.~Long, R.~Tempone, Fast {B}ayesian
  experimental design: Laplace-based importance sampling for the expected
  information gain, Computer Methods in Applied Mechanics and Engineering 334
  (2018) 523--553.

\bibitem{menzies2016efficient}
N.~A. Menzies, An efficient estimator for the expected value of sample
  information, Medical Decision Making 36~(3) (2016) 308--320.

\bibitem{jalal2018gaussian}
H.~Jalal, F.~Alarid-Escudero, A {G}aussian approximation approach for value of
  information analysis, Medical Decision Making 38~(2) (2018) 174--188.

\bibitem{heath2018efficient}
A.~Heath, I.~Manolopoulou, G.~Baio, Efficient {M}onte {C}arlo estimation of the
  expected value of sample information using moment matching, Medical Decision
  Making 38~(2) (2018) 163--173.

\bibitem{heath2020calculating}
A.~Heath, N.~Kunst, C.~Jackson, M.~Strong, F.~Alarid-Escudero, J.~D.
  Goldhaber-Fiebert, G.~Baio, N.~A. Menzies, H.~Jalal, Calculating the expected
  value of sample information in practice: considerations from 3 case studies,
  Medical Decision Making 40~(3) (2020) 314--326.

\bibitem{duffie2010dynamic}
D.~Duffie, Dynamic asset pricing theory, Princeton University Press, 2010.

\bibitem{gordy2010nested}
M.~B. Gordy, S.~Juneja, Nested simulation in portfolio risk measurement,
  Management Science 56~(10) (2010) 1833--1848.

\bibitem{broadie2015risk}
M.~Broadie, Y.~Du, C.~C. Moallemi, Risk estimation via regression, Operations
  Research 63~(5) (2015) 1077--1097.

\bibitem{hong2017kernel}
L.~J. Hong, S.~Juneja, G.~Liu, Kernel smoothing for nested estimation with
  application to portfolio risk measurement, Operations Research 65~(3) (2017)
  657--673.

\bibitem{giles2019multilevel}
M.~B. Giles, A.-L. Haji-Ali, Multilevel nested simulation for efficient risk
  estimation, SIAM/ASA Journal on Uncertainty Quantification 7~(2) (2019)
  497--525.

\end{thebibliography}

\end{document}